%% file: dcorr4.tex
\documentclass[12pt]{article}
\usepackage{amsmath,amsthm,amssymb}
\newtheorem{theorem}{Theorem}
\newtheorem{lemma}[theorem]{Lemma}

\begin{document}

\title {\large \sc Neglecting Discretization Corrections in\\
         Regularized Singular or Nearly Singular Integrals
 }

\author{
{\normalsize J. Thomas Beale} \\
{\normalsize{\em Department of Mathematics, Duke University, Box 90320}}\\
{\normalsize{\em Durham, North Carolina 27708, U.S.A.}}\\
}

\date{}

\maketitle
\input deffs.tex

\begin{abstract}
A method for computing singular or nearly singular integrals on closed surfaces
was presented by
J. T. Beale, W. Ying, and J. R. Wilson [Comm. Comput. Phys. 20 (2016), 733--753,  arxiv.org/abs/1508.00265] and applied to single and double layer potentials for
harmonic functions.  It uses regularized kernels, a straightforward quadrature rule,
and corrections added for smoothing and discretization errors.  In this note we
give estimates for the discretization corrections which show that they can reasonably be neglected with proper choice of numerical parameters.
\end{abstract}

A method for computing singular or nearly singular integrals is presented in [1] and applied to single and double layer potentials for harmonic functions.  The potential evaluated on the surface is a singular integral, and when evaluated off the surface but nearby, it is nearly singular.  The procedure has these steps: (1) replace the singular kernel with a regularized or smooth version; (2) compute a sum determined by a quadrature rule for surface integrals; and (3) add two corrections, one for the smoothing error and one for the discretization error.
In this note we show that it is reasonable to neglect the discretization correction, with proper choice of the numerical parameters.  Thus the amount of work needed is reduced. Even though the correction is first order in the mesh spacing $h$, the coefficients are quite small.  Of course we cannot say the entire discretization error is small without qualification.  The remaining error is higher order in $h$ but can be large, e.g. if the curvature of the surface is large.

The numerical parameters are the mesh spacing $h$; the smoothing radius $\delta$, or equivalently the ratio $\rho = \delta/h$; the angle $\theta$ in the partition of unity for the unit sphere used in the surface quadrature rule; and a coefficient $a$ in the partition of unity explained below.  In [1] we used $\rho = 2$ for the nearly singular case off the surface and $\rho = 3$ for the method on the surface.  We used $\theta = 70^o$ and $a = 1$.  Here we recommend the same choices, except that we find the choice $a = 2$ can lead to smaller coefficients in the discretization correction than $a = 1$.  We do not recommend $a > 2$ because of larger derivatives in the partition of unity.

Assuming we use the procedure of [1] except that we neglect the discretization correction, the remaining error in the nearly singular case has the form
$ c_1 h + c_2 h^2 + c_3\delta^3 $.  Here $c_1$ and $c_2$ depend on $\rho$.
The first term comes from the neglected correction.  
In this note we give specific estimates for $c_1$.
In [1] we remarked that $c_2$ is small, depending of course on the smoothness of the underlying problem.  Since $c_1$ and $c_2$ are small, in practical terms we see
errors about $O(h^3)$ as we decrease $h$ while holding $\rho = \delta/h$ constant,
e.g. $\rho = 2$,
but that cannot be true as $h \to 0$.  Alternatively we could choose
$\delta = \delta_0 h^q$ with $q < 1$, so that $\rho$ increases as $h$ decreases.
The error will then be $O(\delta^3) = O(h^{3q})$ as $h \to 0$, because $c_1$ and $c_2$
contain gaussian factors which decrease rapidly as $\rho$ increases.  For
evaluation {\it on} the surface, the exactly singular case, we use special kernels, so that 
the method has a higher order smoothing error,
and the error without any corrections has the form $c_1h + c_2h^2 + c_3\delta^5$.
The same remarks as before apply to the low order terms $c_1 h + c_2 h^2$.
With $\rho = 3$ we see smaller errors than in the nearly singular case.  
In principle we can again set $\delta = \delta_0 h^q$ with $q<1$ and obtain $O(h^{5q})$ convergence.
The estimates described here were derived in [2].

In [4] S. Tlupova and the author used a related approach to compute single or double layer integrals for Stokes flow.  We found the expected order of convergence in computational examples even though discretization corrections were not included.  The results in this note do not apply to the Stokes case but do support the expectation that the discretization corrections should not be needed.

\medskip

{\bf Statement of results.}
We first summarize the conclusions and then give details.
Equation numbers here refer to [1].
For the single layer integral (3.1) off the surface but nearby, with density function $\psi$, the discretization correction ${\cal T}_2$ is given by (3.6).  We show here that
\beq  |{\cal T}_2| \lee \eps_0 h \max|\psi|  \eeq
where $\eps_0$ is a small number depending on $\rho$ and $a$.  Values of $\eps_0$ are given in Table 1. Here $\theta = 70^o$.  As usual $e\!\!-\!\!d$ means $10^{-d}$.
 
\begin{table}[ht]
\caption{Coefficient $\eps_0$ for the single layer correction (3.6).} 
\vspace{4pt}
\centering 
\begin{tabular}{ c | c c c } 
 & $\rho = 2$ & $\rho = 2.5$ & $\rho = 3$  \\ \hline 
  $a = 1$ & 1.8e-6 & 1.6e-8 & 7.7e-11 \\
  $a = 2$ & 2.1e-7 & 6.0e-10 & 1.2e-12 
\end{tabular}
\end{table}

For the single layer potential evaluated {\it on} the surface, we use a more special regularized kernel with higher order smoothing error and $\rho = 3$.  The discretization correction is (3.15).  It has an estimate with the same qualitative form (1), and we find
\beq \eps_0 = 9.3e\!\!-\!\!8\;\;\mbox{for}\;\;a = 1 \,,\quad 
   \eps_0 = 1.9e\!\!-\!\!9\;\;\mbox{for}\;\;a = 2  \eeq

For the double layer integral (3.7) off the surface the discretization correction
${\cal N}_2$ is given by (3.12),(3.13).  We show that
\beq  |{\cal N}_2| \lee \eps_1 h \max|\nabla\varphi|  \eeq
where $\nabla\varphi$ is the tangential gradient of the density function $\varphi$ and
$\eps_1$ depends on $\rho$ and $a$.  Values of $\eps_1$ are given in Table 2, again with $\theta = 70^o$. 

\begin{table}[ht]
\caption{Coefficient $\eps_1$ for the double layer correction (3.12),(3.13).} 
\vspace{4pt}
\centering 
\begin{tabular}{ c | c c c } 
 & $\rho = 2$ & $\rho = 2.5$ & $\rho = 3$  \\ \hline 
   $a = 1$ & 7.0e-6 & 7.6e-8 & 4.6e-10 \\
   $a = 2$ & 8.3e-7 & 3.0e-9 & 6.6e-12
\end{tabular}
\end{table}

For the double layer potential {\it on} the surface, no discretization correction was needed; see Sec. 3.3 of [1].


\medskip

{\bf Properties of the function E.}
The corrections use the function
\beq E(p,q) = e^{2pq}\text{erfc}(p + q) + e^{-2pq}\text{erfc}(-p + q) \eeq
where $\text{erfc}$ is the complementary error function.  We collect some properties in a lemma, proved below.

\begin{lemma}
(1) For all $z$, $|\mbox{erfc}(z)| \leq 2$, and for $z \geq 0$,
\beq \mbox{erfc}(z) \lee e^{-z^2}/(1 + z)\,,\quad z \geq 0 \eeq.

(2)  $E$ is even in $p$ and positive.  For $p \geq 0$ and $q > 0$ it is decreasing in
either $p$ or $q$. 

(3) Assume $p \geq 0$ and $q > 0$.  Then
\beq E(p,q) \lee \left((1+q)^{-1} + 1\right) e^{-p^2} e^{-q^2} 
              \,,\quad  p \leq q \eeq
\beq E(p,q) \lee 3e^{-pq}e^{-q^2}   \,,\quad  p \geq q \eeq

(4) If $p \geq q > 1$, $pE$ is decreasing in $p$.

(5) For each $q>1$ and all $p \geq 0$,
\beq pE(p,q) \leq .429 \left((1+q)^{-1} + 1\right) e^{-q^2} \eeq
\end{lemma}

{\bf Role of the partition of unity on the sphere.}
The corrections have factors with the partition of unity functions $\zeta^k$.
The $\zeta^k$ are functions of the normal vector, with
$\zeta^k = \beta_k/\beta_{tot}$ where $\beta_k$ depends on $n_k$,
the $k$th component of the normal,
and $\beta_{tot} = \beta_1 + \beta_2 + \beta_3$.
With an angle $\theta$ and parameter $a$ chosen, 
\beq \beta_k = \exp(ar^2/(r^2-1))\,, \quad r = (\cos^{-1}|n_k|)/\theta \eeq
provided $|n_k| > \cos{\theta}$, and zero otherwise.
In [1], $a = 1$, but here we allow other choices of $a$.
For the function $n_k \mapsto \beta_k$ with $\theta = 70^o$, we find the maximum
derivative is $2.3$ with $a = 1$, $2.7$ with $a = 2$, and larger for $a>2$.
For this reason we recommend using $a \leq 2$.

The corrections have sums over integer pairs $(m_1,m_2)$ in 
\beq  Q \,\equiv\, \{(m_1,m_2)\in \Z^2: m_2 > 0\; \mbox{or}\;
     (m_2=0\; \mbox{and} \; m_1>0)\} \eeq
(In [1] we used $n$ rather than $m$; it is changed here to avoid confusion with the normal vector.)  The corrections use
$\|m\|^2$, defined as $\Sigma_{i,j} g^{ij}m_im_j$. 
Here $g^{ij}$ is the inverse metric tensor in the $k$th coordinate patch. 
From formulas in the Monge patch $z = f(x,y)$ we see that $\|m\| \geq \gamma_3 |m|$, where
$|m|$ is the Euclidean norm of $m = (m_1,m_2)$
and $\gamma_3 = (1 + f_x^2 + f_y^2)^{-1/2}$,
which is the same as the absolute value of the third component of the normal vector; similarly for the other two cases.  In summary, for the $k$th patch, we have
$\|m\| \geq \gamma_k |m|$ with $\gamma_k = |n_k|$.  We can assume $\gamma_k \geq
\gamma_0 \equiv \cos{\theta}$, since $\beta_k = 0$ for $|n_k| \leq \cos{\theta}$.
The correction formulas have factors $\beta_k$ times $E(p,q)$ with
$q \geq \pi\rho\gamma_k|m|$.  A key point is that $E$ might be relatively large if
$\gamma_k$ is small, but then $\beta_k$ is small.  Thus the product is smaller than it might appear at first.  This was not evident in [1].

\medskip  

{\bf The single layer correction.}
We now discuss the correction (3.6) for the single layer off the surface.  Because of Lemma 1(2) and the remarks above, we can majorize $E(\lambda,\pi\rho\|m\|)$ by
$E(0,\pi\rho\gamma_k|m|)$ and the sum multiplying $h\psi$ by
\beq  \frac{1}{4\pi}  \sum_{m \in Q} \sum_{k=1}^3
      \frac{\zeta_k}{\gamma_k |m|}E(0,\pi\rho\gamma_k|m|)  \eeq
The sum depends on the normal vector or equivalently on the unit vector
$\gamma_k$, $k=1,2,3$.  To find the upper bounds in Table 1 we truncate the sum
(e.g., $m_1,m_2 \leq 2$), compute directly for various $\gamma_k$, and maximize over admissible choices of $\gamma_k$.  The terms $m = (1,0)$ and $(0,1)$ are much larger than the others.  For remaining $m$ we use 
\beq  \gamma_k \geq \gamma_0 \equiv \cos{\theta}\,,
                  \qquad q_0 \equiv \pi\rho\gamma_0  \eeq
and Lemma 1(2) to estimate the $m$th term above by
        $(4\pi\gamma_0 |m|)^{-1} E(0,q_0|m|)$.
Finally, we use this crude estimate for the terms with higher $m$ to show their sum is negligible.
For this last step we choose $R$, e.g., $R = 2$, and define $Q_R$ as the subset of $Q$ with $(m_1,m_2)$ such that $m_1,m_2 > 0$ and $|(m_1-1,m_2-1)| \geq R$;
or $m_1-1 \geq R$ and $m_2 = 0$; or $m_2-1 \geq R$ and $m_1 = 0$.  We bound the sum
over $Q_R$ with an integral and obtain the following estimate.

\begin{lemma}The 
contribution to the part of the sum in (11) from $m \in Q_R$ is bounded by
\beq \frac{\sqrt{\pi}}{\gamma_0 q_0^2}\text{erfc}(q_0 R)
        \left( \frac{1}{4R} + \frac{1}{2\pi R^2}\right)  \eeq
\end{lemma}

For the case of evaluation {\it on} the surface, the correction (3.15) is treated similarly.  Again the main contribution is from the terms $m = (1,0)$ and
$(0,1)$, and we check that the rest is negligible.  In (3.15) the part with erfc is the same as before.  The second part gives an additional term in the estimate for large $m$:

\begin{lemma}
For (3.15) the sum over $Q_R$, as in Lemma 2, is bounded by the previous term plus
\beq \frac{\rho^2}{2\sqrt{\pi}q_0^2}\left(\pi + \frac{2}{R}\right)
  \left(\frac53 + \frac23 q_0^2R^2\right)\,e^{-q_0^2R^2}  \eeq
\end{lemma}

\medskip

{\bf The double layer correction.}
We handle the correction (3.12), (3.13) for the double layer in a similar way, except that we cannot easily identify the maximum in $\lambda$. The correction includes a sum
$\Sigma_{r,s = 1}^2 (\pa_r\varphi)g^{rs}m_s$  where $g^{rs}$ is the inverse metric tensor in one of the three coordinate systems and $\pa_r$, $r=1,2$ are the partial derivatives in the coordinate system.  The tangential gradient
$\nabla\varphi$ of $\varphi$ is
$\Sigma_{r,s} (\pa_r\varphi) g^{rs} T_s$; it is independent of coordinates.  Here
$T_s$ is the $s$th tangent vector in the $k$th system.  If e.g. $k=3$, the surface
has the form $x_3 = z(x_1,x_2)$, and $T_1 = (1,0,z_1)$, $T_2 = (0,1,z_2)$.  In this case the $s$-component of $\nabla\varphi$ is $\Sigma_{r,s}(\pa_r\varphi) g^{rs}$,
$s = 1,2$, and thus $\Sigma_{r,s = 1}^2 (\pa_r\varphi)g^{rs}m_s =
\Sigma_{s=1}^2 (\nabla\varphi)_s m_s$.  In particular 
$|\Sigma_{r,s = 1}^2 (\pa_r\varphi)g^{rs}m_s| \leq |\nabla\varphi||m|$.
Again we use $\|m\| \geq \gamma_k|m| \geq \gamma_0|m|$.

We treat the terms $m = (1,0)$ and $(0,1)$ in (3.12), (3.13) more carefully than the others, since they are again the largest.  We take absolute values and use
$|(\nabla\varphi)_1| + |(\nabla\varphi)_2| \leq \sqrt{2} M_1$, where
$M_1 = \max{|\nabla\varphi|}$.
The sum of these two terms is bounded by
$hM_1$ times
\beq \frac{\sqrt{2}\rho\lambda}{2} \sum_{k=1}^3
      \frac{\zeta_k}{\gamma_k}E(\lambda,\pi\rho\gamma_k)  \eeq
We find an upper bound for this quantity by
direct computation, maximizing over admissible $\gamma_k$ and
$\lambda$ with $0 \leq \lambda \leq q_0$;
we know from Lemma 1(4) that the maximum occurs for $\lambda$ in this interval.
The maxima occur for $\lambda$ between $.7$ and $.8$.      

It then remains to check that the other terms are negligible in comparison.
For other $m$ we have $|\Sigma_{r,s}\pa\varphi_r g^{rs}m_s|/|m| \leq M_1$ and
from Lemma 1(5), $\lambda E(\lambda,q_0|m|) \leq .86 e^{-q_0^2|m|^2}$.
Thus the $m$th term is bounded by $hM_1$ times 
\beq (\rho/2\gamma) e^{-q_0^2|m|^2}  \eeq 
We verify computationally that a few terms are negligible, and we can estimate the remaining ones as in Lemma 2, comparing the sum of these terms with an integral:

\begin{lemma}
The sum over $Q_R$ of the terms (16) is bounded by
\beq \frac{\rho}{4\gamma q_0^2}(\pi + \frac{2}{R})e^{-q_0^2R^2}  \eeq
\end{lemma}

\medskip

{\bf Proof of Lemma 1.}
(1)  The inequality is 
derived from the integral formula (7.7.1) for $\text{erfc}$ in the NIST
DLMF,
\beq \text{erfc}(z) \equiv \frac{2}{\sqrt{\pi}}\int_z^\infty e^{-t^2}\,dt
       \eq \frac{2}{\pi}e^{-z^2}
            \int_0^\infty \frac{e^{-z^2t^2} }{t^2 + 1}\,dt  \eeq
using the inequality $\exp(-z^2t^2) \leq 1/(z^2t^2 + 1)$, which
leads to a contour integral.

(2) The first statement is obvious from the definition.
We check that $\pa E/\pa q < 0$
assuming $p \geq 0$ and $q > 0$.
The partial derivative is
\beq 2pe^{2pq}\text{erfc}(p+q) - 2pe^{-2pq}\text{erfc}(-p+q)
 -(2/\sqrt{\pi})e^{2pq}e^{-(p+q)^2} - (2/\sqrt{\pi})e^{-2pq}e^{-(-p+q)^2} \eeq
This has the form $T_1 - T_2 - T_3 - T_4$ with each $T_i \geq 0$.
We check that $T_1 < T_3 + T_4$ so that the sum is $<0$.  Combining exponents in $T_3$ and $T_4$ we see they are the same, and
$ T_3 + T_4 > 2e^{-p^2-q^2} $ since $2/\sqrt{\pi} > 1$.  For $T_1$ we use
the estimate for $\text{erfc}$ to get 
\beq T_1 \leq 2pe^{2pq}(1+p+q)^{-1}e^{-(p+q)^2} \leq 2p(p+1)^{-1}e^{-p^2-q^2} \eeq
verifying that $T_1 < T_3 + T_4$.

Next we check that $\pa E/\pa p < 0$.  We find two terms cancel in this derivative, and
\beq \frac{\pa E(p,q)}{\pa p} \eq 2qE^{-}(p,q) \eeq
with 
\beq E^{-}(p,q) = e^{2pq}\text{erfc}(p + q) - e^{-2pq}\text{erfc}(-p + q) \eeq
Computing $\pa E^{-}/\pa q$ as for $E$, we again find that two terms cancel,
and $\pa E^{-}/\pa q = 2pE > 0$, so that $E^{-}$ increases in $q$.  Since $E^{-} \to 0$ as $q \to \infty$, we see that $E^{-} < 0$ and also
$\pa E/\pa p < 0$.

(3) With $p,q \geq 0$, the
first term in $E$ is always bounded by $\exp(-p^2-q^2)/(1 + q)$.  If
$p \leq q$, the second term is bounded by $\exp(-p^2-q^2)$.  Thus
\beq E(p,q) \lee \left((1+q)^{-1} + 1\right) e^{-p^2} e^{-q^2} 
              \,,\quad  p \leq q \eeq
Suppose instead that $p \geq q$.  Then the first
term in $E$ is bounded by $\exp(-p^2-q^2) \leq \exp(-pq)\exp(-q^2)$. 
The second term is bounded by $2\exp(-2pq) = 2\exp(-pq)\exp(-pq)
\leq 2 \exp(-pq)\exp(-q^2)$.  Thus 
\beq E(p,q) \lee 3e^{-pq}e^{-q^2} \,, \quad p \geq q  \eeq

(4) We write
$E = E_1 + E_2$ and $E^{-} = E_1 - E_2$.  Then 
$\pa(pE)/\pa p = E + 2pqE^{-} = E_1 + E_2 + 2pqE_1 - 2pqE_2$.
We need to show that $\mu E_2 > E_1$ with $\mu = (2pq-1)/(2pq+1)$.
We start by noting as above that $E_1 \leq \nu e^{-(p^2+q^2)}$ with
$\nu = 1/(1 + p + q)$.  Since $p \geq q$, $E_2 \geq e^{-2pq}\cdot 1$.  
Now we need to show $\mu e^{-2pq} \geq \nu e^{-(p^2+q^2)}$.  Since
$ 2pq \leq p^2 + q^2$, $e^{-2pq} \geq e^{-(p^2+q^2)}$.  It remains to show
that $\mu \geq \nu$.  For $p = q = 1$, $\mu = \nu$, and as either $p$ or $q$ increases beyond $1$, $\mu$ increases and $\nu$ decreases, so this last inequality holds for $p,q \geq 1$.

(5) For fixed $q$ the estimate for $pE$ with $p \leq q$,
using (23) above for $E$, is maximized
at $p = 1/\sqrt{2}$, and
\beq pE(p,q) \leq .429(1 + (1+q)^{-1})e^{-q^2} \eeq
We know from (4) that for $p \geq q > 1$, $pE$ decreases in $p$, so this inequality holds
for all $p$.

\medskip

{\bf Proof of Lemma 2.}  We need an upper bound for the sum over $m \in Q_R$ of
\beq (4\pi\gamma_0 |m|)^{-1} E(0,q_0|m|) = 
           (2\pi\gamma_0 |m|)^{-1} \text{erfc}(q_0|m|)  \eeq
The term is a decreasing function of
$|m|$.  We bound the sum with $m_1,m2>0$ by a double integral and the
(equal) sums over $(m_1,0)$ and $(m_2,0)$ by single integrals; the sum is bounded by
$(I_1 + 2I_2)/(2\pi\gamma_0)$, where
\beq I_1 = \pi \int_R^\infty r^{-1} \text{erfc}(q_0r) r\,dr\,, \quad
    I_2 = \int_R^\infty r^{-1} \text{erfc}(q_0r)\,dr  \eeq
For $I_1$ we use $\text{erfc}(q_0|m|) \leq \exp(-q_0^2|m|^2)/(q_0|m|)
 \leq \exp(-q_0^2|m|^2)/(q_0 R)$, change to 
$s = q_0 r$, and integrate to find
$I_1 \leq (\pi/(q_0^2R))(\sqrt{\pi}/2)\text{erfc}(q_0R)$.
Inside the integral $I_2$ we use $1/r \leq 1/R$ leading to the same integral as $I_1$.

\bigskip

{\bf References}

\smallskip

[1] J. T. Beale, W. Ying, and J. R. Wilson, {\em A simple method for computing singular or nearly singular integrals on closed surfaces},  Comm. Comput. Phys. 20 (2016), 733--753 or  arxiv.org/abs/1508.00265.

[2] J. T. Beale, {\em A grid-based boundary integral method for elliptic problems in three dimensions}, SIAM J. Numer. Anal. 42 (2004), 599--620.
 
[3] J. T. Beale and M. C. Lai, {\em A method for computing nearly singular integrals},
SIAM J. Numer. Anal. 38 (2001), 1902--1925. 

[4] S. Tlupova and J. T. Beale, {\em Regularized single and double layer integrals in 3D Stokes flow},  J. Comput. Phys.  386 (2019), 568--584 or arxiv.org/abs/1808.02177.

\end{document}

%% file: deffs.tex
\newcommand{\beq}{\begin{equation}}
\newcommand{\eeq}{\end{equation}}
\newcommand\p{\,+\,}
\newcommand\lee{\,\leq\,}
\newcommand\gee{\,\geq\,}
\newcommand\m{\,-\,}
\newcommand\eq{\,=\,}
\newcommand{\eps}{\varepsilon}
\newcommand{\sig}{\sigma}
\newcommand{\pa}{\partial}
\newcommand{\lilhalf}{{\textstyle \frac12}}
\newcommand{\lilth}{{\textstyle \frac32}}
\newcommand{\mm} {_{\max}}
\newcommand{\np}{{n+1}}
\newcommand{\nm}{{n-1}}
\newcommand{\nhalf}{{n+1/2}}
\newcommand{\ex}{{ex}}
\newcommand{\NN}{{\mathcal N}}
\newcommand{\R}{{\mathbb R}}
\newcommand{\Z}{{\mathbb Z}}
\newcommand{\laph}{\Delta_h}
\newcommand{\Ptw}{{\tilde P}}
\newcommand{\fcns}{{\mathcal F}(\Omega_h)}
\newcommand{\phihat}{{\hat \varphi}}
\newcommand{\fhat}{{\hat f}}
\newcommand{\ghat}{{\hat g}}
\newcommand{\ip}{{i+1}}
\newcommand{\im}{{i-1}}
\newcommand{\iph}{{i+1/2}}
\newcommand{\imh}{{i-1/2}}
\newcommand{\jp}{{j+1}}
\newcommand{\jm}{{j-1}}
\newcommand{\jph}{{j+1/2}}
\newcommand{\jmh}{{j-1/2}}

\newcommand{\kp}{{k+1}}
\newcommand{\km}{{k-1}}
\newcommand{\kph}{{k+1/2}}
\newcommand{\kmh}{{k-1/2}}
\newcommand{\om}{\Omega}

\newcommand{\xtw}{\tilde x}
\newcommand{\ytw}{\tilde y}
\newcommand{\utw}{\tilde u}
\newcommand{\vtw}{\tilde v}
\newcommand{\wtw}{\tilde w}
\newcommand{\phitw}{\tilde \varphi}

\newcommand{\bbox}{\cal B}
\newcommand{\boxh}{{\cal B}_h}
\newcommand{\omp}{\om^+}
\newcommand{\omm}{\om^-}
\newcommand{\omhp}{\om_h^+}
\newcommand{\omhm}{\om_h^-}
\newcommand{\omhpm}{\om_h^\pm}
\newcommand{\ompcl}{\overline \omp}
\newcommand{\ommcl}{\overline \omm}
\newcommand{\omhpcl}{\overline \omhp}
\newcommand{\omhmcl}{\overline \omhm}
\newcommand{\omhpmcl}{\overline \omhpm}

\newcommand{\fcnspeq}{{\cal F}_{\#}(\omhpcl)}
\newcommand{\fcnseq}{{\cal F}_{\#}(\Gamma_h^1)}
\newcommand{\fcnscut}{{\cal F}(\Gamma_h^0)}
\newcommand{\fcnsgam}{{\cal F}(\Gamma_h^1)}
\newcommand{\fcnsp}{{\cal F}(\omhpcl)}
\newcommand{\fcnsm}{{\cal F}(\omhmcl)}
\newcommand{\fcnspm}{{\cal F}(\omhpmcl)}
\newcommand{\fcnspz}{{\cal F}_0(\omhpcl)}
\newcommand{\ints}{\cal I}
\newcommand{\intsh}{{\cal I}_h}
\newcommand{\intscut}{{\cal I}^0}
\newcommand{\intshcut}{\intscut_h}
\newcommand{\intsp}{{\ints^+}}
\newcommand{\intsm}{{\ints^-}}

\newcommand{\intshp}{{\cal I}^+_h}
\newcommand{\intshm}{{\cal I}^-_h}

\newcommand{\A}{{\cal A}}
\newcommand{\B}{{\cal B}}

\newcommand{\fbar}{\overline f}